\documentclass[11pt, letterpaper]{amsart}
\usepackage{mathbbol}
\usepackage{amsmath, amsthm, amssymb}
\usepackage{txfonts}

\newtheorem{thm}{Theorem}[section]
\newcommand{\bt}{\begin{thm}}
\newcommand{\et}{\end{thm}}

\newtheorem{cor}[thm]{Corollary}
\newcommand{\bc}{\begin{cor}}
\newcommand{\ec}{\end{cor}}

\newtheorem{lem}[thm]{Lemma}
\newcommand{\bl}{\begin{lem}}
\newcommand{\el}{\end{lem}}

\newtheorem{prop}[thm]{Proposition}
\newcommand{\bp}{\begin{prop}}
\newcommand{\ep}{\end{prop}}

\newtheorem{defn}[thm]{Definition}
\newcommand{\bd}{\begin{defn}}      
\newcommand{\ed}{\end{defn}}

\newtheorem{rmrk}[thm]{Remark}
\newcommand{\br}{\begin{rmrk}}
\newcommand{\er}{\end{rmrk}}

\newcommand{\thmref}[1]{Theorem~\ref{#1}}
\newcommand{\secref}[1]{Section~\ref{#1}}
\newcommand{\lemref}[1]{Lemma~\ref{#1}}
\newcommand{\defref}[1]{Definition~\ref{#1}}

\newcommand{\propref}[1]{Proposition~\ref{#1}}
\newcommand{\remref}[1]{Remark~\ref{#1}}

\newcommand{\N}{\mathbb{N}}

\newcommand{\R}{\mathbb{R}}
\newcommand{\E}{\mathbb{E}}
\newcommand{\Z}{\mathbb{Z}}

\newcommand{\dist}{\operatorname{dist}}
\newcommand{\diam}{\operatorname{diam}}
\newcommand{\length}{\operatorname{length}}
\newcommand{\area}{\operatorname{Area}}
\newcommand{\hm}{{\mathcal H}}
\newcommand{\lm}{{\mathcal L}}
\newcommand{\md}{\operatorname{md}}
\newcommand{\jac}{{\mathbf J}}

\newcommand{\image}{\operatorname{Im}}

\newcommand{\lip}{\operatorname{Lip}}
\newcommand{\mass}[2][]{{\mathbf M_{#1}}(#2)}

\newcommand{\intcurr}{{\mathbf I}}      

\newcommand{\fillarea}{{\operatorname{Fill\,Area}}}
\newcommand{\fillrad}{{\operatorname{Fill\,Rad}}}

\newcommand{\rstr}{\:\mbox{\rule{0.1ex}{1.2ex}\rule{1.1ex}{0.1ex}}\:}
\newcommand{\bdry}{\partial}

\newcommand{\on}[1]{|_{#1}}
\newcommand{\spt}{\operatorname{spt}}
\newcommand{\ohne}{\backslash}

\begin{document}

\title{Gromov hyperbolic spaces and the sharp isoperimetric constant}

\author{Stefan Wenger}

\address
  {Courant Institute of Mathematical Sciences\\
   251 Mercer Street\\
   New York, NY 10012}
\email{wenger@cims.nyu.edu}

\date{March 19, 2007}


\begin{abstract}
 In this article we exhibit the largest constant in a quadratic isoperimetric inequality 
 which ensures that a geodesic metric space is Gromov hyperbolic. As a particular consequence we obtain that
 Euclidean space is a borderline case for Gromov hyperbolicity in terms of the isoperimetric function.  We prove similar results for the linear filling radius inequality. 
 Our results strengthen and generalize theorems of Gromov, Papasoglu and others.
 \end{abstract}

\maketitle

\bigskip

\section{Introduction}
The classical isoperimetric inequality in the Euclidean plane $\E^2$ asserts that the area $A$ enclosed by a closed curve $\gamma$ in $\E^2$ satisfies
\begin{equation*}
 A\leq \frac{1}{4\pi}\length(\gamma)^2,
\end{equation*}
with equality if and only if $\gamma$ parametrizes a circle. One of the main purposes of the present article is to prove the sharp result below, which shows
that a geodesic metric space cannot have a quadratic isoperimetric inequality with constant strictly smaller than $\frac{1}{4\pi}$
unless it is Gromov hyperbolic (and thus already admits a (coarse) linear isoperimetric inequality).


\bt\label{theorem:isop-special-case}
 Let $X$ be a geodesic metric space and suppose there exists $\varepsilon>0$ such that every sufficiently long Lipschitz loop $\gamma$ in $X$
 bounds a singular Lipschitz disc $\Sigma$ in $X$ of area
 \begin{equation}\label{equation:isop-special-case}
  \area(\Sigma)\leq \frac{1-\varepsilon}{4\pi}\length(\gamma)^2.
 \end{equation}
 Then $X$ is Gromov hyperbolic.
\et
More general results will be described below and in \secref{section:isop-gromov}.
By definition, a singular Lipschitz disc in $X$ is (the image of) a Lipschitz map $\varphi: D^2\to X$, where $D^2\subset\E^2$ is the unit disc. 
Furthermore, $\area(\Sigma)$ is the `parametrized' $2$-dimensional Hausdorff measure of $\Sigma$, see \secref{section:metric-derivatives}. In particular, if $\varphi$ is one-to-one on a set of full measure then $\area(\Sigma)=\hm^2(\Sigma)$, where $\hm^2$ is the $2$-dimensional Hausdorff measure on $X$.

Recall that, by definition, a geodesic metric space $X$ is $\delta$-hyperbolic if every geodesic triangle in $X$ is 
$\delta$-slim, i.e.\ if each side of the triangle is contained in the $\delta$-neighborhood of the union of the other two sides. 
The theory of $\delta$-hyperbolic spaces (and groups) goes back to Gromov \cite{Gromov-hyperbolic}. A geodesic metric space is therefore said to be Gromov hyperbolic if it is 
$\delta$-hyperbolic for some $\delta\geq0$. 
It is well-known that Gromov hyperbolic spaces admit a (coarse) linear isoperimetric inequality for curves. More precisely,
if $X$ is $\delta$-hyperbolic and if there exists $C>0$ such that every Lipschitz loop in $X$ of length at most $20\delta$ bounds a singular Lipschitz disc of area at most
$C$, then $X$ admits a linear isoperimetric inequality for curves, i.e.\ every Lipschitz loop $\gamma$ in $X$ bounds a singular Lipschitz disc $\Sigma$ with
\begin{equation*}
 \area(\Sigma)\leq D\length(\gamma),
\end{equation*}
where $D$ only depends on $C$ and $\delta$.

Clearly, the constant $\frac{1}{4\pi}$ appearing in \eqref{equation:isop-special-case} is optimal as follows from the classical isoperimetric inequality in $\E^2$.
\thmref{theorem:isop-special-case} is new even in the setting of Riemannian manifolds and was previously only known in the special case when $X$ is a Hadamard manifold or, 
more generally, a ${\rm CAT}(0)$-space (for which it was observed by Gromov). 
In the special setting of Riemannian manifolds the best  constant previously established was $\frac{1}{16\pi}$, due to Gromov \cite{Gromov-hyperbolic}.
Indeed, using conformal mappings Gromov proved that a `reasonable' Riemannian manifold $M$ is $\delta$-hyperbolic provided \eqref{equation:isop-special-case} holds 
with $\varepsilon:= \frac{3}{4}$, i.e.\ if
every sufficiently long Lipschitz loop $\gamma$ in $M$ bounds a singular Lipschitz disc $\Sigma$ in $M$ of area
\begin{equation*}
 \area(\Sigma)\leq \frac{1}{16\pi}\length(\gamma)^2.
\end{equation*}
For the meaning of `reasonable' see \cite[p.\ 176]{Gromov-hyperbolic}. For example, the universal covering of a closed Riemannian manifold is `reasonable'. 
See also \cite{Coornaert-Delzant-Papadopoulos}, where a detailed account of Gromov's proof is given.
Gromov furthermore showed that the same conclusion holds for geodesic metric spaces provided \eqref{equation:isop-special-case} is satisfied with $\varepsilon\in(0,1)$ close enough to
$1$. Similar results and alternative proofs of the latter were later given by Olshanskii 
\cite{Olshanskii}, Short \cite{Short}, Bowditch \cite{Bowditch-subquadratic},  Papasoglu \cite{Papasoglu}, and Dru\c{t}u \cite{Drutu-french}. We refer to \cite{Drutu} for an account of the existing results.

In actuality, \thmref{theorem:isop-special-case} is merely a special case of the main result of this paper, \thmref{theorem:isop-general-case}, which will be given in 
\secref{section:isop-gromov}.
To give a rough description of the main theorem let $X$ be a geodesic metric space and $\gamma$ a Lipschitz loop in $X$. Given a metric space
$Y$ in which $X$ isometrically embeds, the filling area of $\gamma$ in $Y$ is, by definition, the least area of a singular Lipschitz chain in $Y$ with boundary $\gamma$.
Recall that $X$ embeds isometrically into $L^\infty(X)$.
It is not difficult to show, see \lemref{lemma:Linfty-fillings}, that the filling area of $\gamma$ in $L^\infty(X)$ is smaller or equal to that in $Y$ for any $Y$ in which $X$ isometrically 
embeds.
Moreover, since $L^\infty(X)$ is a Banach space, the filling area in $L^\infty(X)$ of $\gamma$ is bounded above by $C\length(\gamma)^2$ for some universal constant $C$, and this
holds even if $\gamma$ does not bound a chain in $X$ (and thus has infinite filling area in $X$).
Our main result then shows that for a large class of geodesic metric spaces (which we will call `admissible') the conclusion of
\thmref{theorem:isop-special-case} holds under the weaker assumption that every sufficiently long Lipschitz loop $\gamma$ in $X$ bounds a singular Lipschitz chain $\Sigma$ in
 $L^\infty(X)$ which satisfies \eqref{equation:isop-special-case}. Note that \thmref{theorem:isop-special-case} asked for a $\Sigma$ {\it in $X$} which satisfies 
 \eqref{equation:isop-special-case}.
For example, length spaces which admit a coarse homological quadratic isoperimetric inequality for curves are admissible. 
This includes in particular Cayley graphs of finitely presented groups with quadratic Dehn function.

The techniques used to prove \thmref{theorem:isop-general-case} furthermore yield an analogous sharp result in terms of filling radius inequalities, see \thmref{theorem:fillrad-hyperbolic} in  \secref{Section:fillrad},
which strengthens and generalizes theorems of Gromov and Papasoglu.

\subsection{Outline of the main argument}\label{section:outline-proof-special-case}
We give a short outline of the proof of \thmref{theorem:isop-special-case}, which is achieved in three steps and is by contradiction.\\
{\it Step1:} In \secref{section:coarse-isoperimetric} it is shown that a geodesic metric space $X$ as in \thmref{theorem:isop-special-case} possesses a thickening $X_\delta$ which 
admits a quadratic isoperimetric inequality for curves. Spaces admitting such a thickening will be called admissible in the sequel, see \defref{definition:admissible}.\\
{\it Step 2:} In \secref{section:asymptotic-subsets} it is proved that if $X$ is not Gromov hyperbolic
then there exist a sequence of sets $Z_n\subset X_\delta$ and numbers $r_n\nearrow\infty$ such that $(Z_n,r_n^{-1}d_{X_\delta})$ converges in the Gromov-Hausdorff sense
 to a compact metric space
$(Z,d_Z)$ which admits a Lipschitz map $\varphi:K\to (Z,d_Z)$ with $K\subset\R^2$ compact and $\hm^2(\varphi(K))>0$. 
The construction of such subsets relies on the quadratic isoperimetric inequality for $X_\delta$ and uses the theory of integral currents in metric spaces, recently developed by Ambrosio
and Kirchheim. 
Roughly speaking, the $Z_n$ are constructed as supports of suitable $2$-dimensional integral currents which, upon rescaling, converge to some limit $S$ in a suitable metric space. 
The assumption that $X$ is not Gromov hyperbolic can be used to show that $S\not=0$. The closure theorem for integral currents shows that $S$ is an integral current
and is thus `parametrized by biLipschitz pieces'. The desired metric space $Z$ is simply the support of $S$.\\
{\it Step 3:} First, we remark that, since the Hausdorff measure bounds the Holmes-Thompson area $\mu^{ht}$ from above, \eqref{equation:isop-special-case} holds after replacing 
$\area(\Sigma)$ by the Holmes-Thompson area $\mathcal{F}_{\mu^{ht}}(\Sigma)$. 
Next, let $Z$ and $\varphi$ be as in step 2. By a Rademacher type differentiability theorem of Kirchheim and 
Korevaar-Schoen it follows that 
$Z$ receives an $(1+\varepsilon')$-biLipschitz copy (with $\varepsilon'>0$ very small) of a piece of a $2$-dimensional normed space $V$, see \thmref{theorem:kirchheim-bilip-finite-set}. 
This is used in \secref{section:isop-gromov}  to construct a closed Lipschitz loop $\gamma$
in $X$ which is $(1+\varepsilon')$-biLipschitz equivalent to the boundary of an isoperimetric subset $\mathbb{I}_V$ of  $V$, i.e.\ a compact convex subset of 
maximal $\mu^{ht}$-measure among all convex subsets with the same perimeter. View $V$ as a linear subspace of $\ell^\infty$. Since $\ell^\infty$ is an injective metric space,
a filling $\Sigma$ satisfying \eqref{equation:isop-special-case} can then be mapped via a $(1+\varepsilon')$-Lipschitz map to a filling of $\bdry\mathbb I_V$ in $\ell^\infty$ which still has
`small' area. 
Since $\mathcal{F}_{\mu^{ht}}$ is semi-elliptic in the class of singular Lipschitz discs (by a recent important result of Burago and Ivanov) and since
\begin{equation*}
 \mu^{ht}(\mathbb{I}_V)= \frac{1}{4\pi}\length(\bdry\mathbb{I}_V)^2
\end{equation*}
this can be shown to lead to a contradiction. 

The change from the Hausdorff measure to the Holmes-Thompson area is made necessary by the fact that it is not known 
whether the Hausdorff measure is semi-elliptic, see \secref{section:area-functionals}. The only reason for using the Hausdorff measure in
\thmref{theorem:isop-special-case} was to make the statement easily accessible.

\subsection{Organization of the paper}
As explained in step 2 of \secref{section:outline-proof-special-case}, the proof of \thmref{theorem:isop-special-case} relies to some extent on the
theory of integral currents in metric spaces. The same applies to Theorems \ref{theorem:isop-general-case} and \ref{theorem:fillrad-hyperbolic}.
For reasons of consistency all results are therefore stated and proved in the language of currents. 
However, apart from those in \secref{section:asymptotic-subsets}, the proofs could be written also in the language of  singular Lipschitz discs or chains (possibly
with some loss in generality). 

Some definitions and facts concerning integral currents in metric spaces are given in \secref{section:preliminaries}. The same section discusses the area
functionals induced by the Hausdorff, the Holmes-Thompson and the Gromov mass$*$ measures.
In \secref{section:coarse-isoperimetric} we show that length spaces which admit a (coarse) quadratic isoperimetric inequality for long curves are admissible in the sense explained
in step 1 of \secref{section:outline-proof-special-case}. As mentioned above, Cayley graphs of finitely presented groups with quadratic Dehn function are admissible.
The principal result of \secref{section:asymptotic-subsets}, described in step 2 of \secref{section:outline-proof-special-case}, is the only place
where results from the theory of metric integral currents enter in a non-trivial way.
The proof of the main result of this article, \thmref{theorem:isop-general-case}, only relies on the results of Sections
\ref{section:asymptotic-subsets} and \ref{section:isop-gromov}, and not on \secref{section:coarse-isoperimetric}.
In \secref{section:isop-gromov} it is furthermore shown how \thmref{theorem:isop-special-case} follows from the main theorem.
Finally, a sharp result involving the filling radius inequality is given in \thmref{theorem:fillrad-hyperbolic}.
 
\medskip

{\bf Acknowledgments:} I would like to thank Juan-Carlos \'Alvarez Paiva, Mario Bonk, Cornelia Dru\c{t}u, Misha Gromov and Bruce Kleiner for discussions and comments. 
Parts of this paper were written during a research visit to the ETH Z\"urich in 2006. I would like to thank the Forschungs\-institut f\"ur Mathematik for its hospitality.

\section{Preliminaries}\label{section:preliminaries}

This section provides definitions and some basic facts concerning $L^\infty$-spaces, injective metric spaces, metric derivatives and
integral currents in metric spaces. The only new result here is \lemref{lemma:almost-optimal-decomposition-one-cycles}.
For background on Gromov hyperbolic spaces we refer to \cite{Gromov-hyperbolic}, \cite{Coornaert-Delzant-Papadopoulos}, \cite{Ghys-delaHarpe}, \cite{Short}.

\subsection{$L^\infty$-spaces, isometric embeddings, and Lipschitz extensions}\label{subsection:injective-spaces}
Given a set $\Omega$ denote by $L^\infty(\Omega)$ the space of bounded $\R$-valued functions on $\Omega$ endowed with the supremum norm
\begin{equation*}
 \|f\|_\infty:= \sup_{a\in\Omega}|f(a)|.
\end{equation*}
We abbreviate $\ell^\infty_n:= L^\infty(\{1,\dots,n\})$ and $\ell^\infty:=L^\infty(\N)$.
As is well-known, if $(Z,d)$ is a metric space and $z_0\in Z$ is fixed, then the map $\varphi_Z: Z\to L^\infty(Z)$ given by $\varphi(z):= d(z,\cdot)-d(z_0, \cdot)$ defines 
an isometric embedding, called the Kuratowski embedding. If $Z$ is separable, then there is an isometric embedding into $\ell^\infty$. 
If $Z$ is a separable Banach space, the isometry may be chosen to be linear, as follows from the Hahn-Banach theorem.

A metric space $X$ is called injective if for every triple $(Z,Y,f)$, where $Z$ is a metric space, $Y$ a subset of $Z$, and $f:Y\to X$ a $1$-Lipschitz
map, there exists an extension $\bar{f}:Z\to X$ of $f$ which is  $1$-Lipschitz. It can be shown (see e.g.~\cite[p.12--13]{Benyamini-Lindenstrauss}) that $L^\infty(\Omega)$ is an
injective space for all sets $\Omega$.
In \cite{Isbell} Isbell furthermore showed that for each metric space $X$ there exists a `minimal' injective space containing $X$. This space is called the 
injective envelope of $X$ and we refer to \cite{Isbell} for its construction and useful properties. 
Injective envelopes will be used in the proof of \propref{proposition:admissible-thickening}.

\subsection{Lipschitz maps and metric derivatives}\label{section:metric-derivatives}
The proof of our main result relies in a crucial way on the following metric differentiability property of Lipschitz maps from Euclidean space into arbitrary metric spaces. 
Let $\varphi: U\to X$ be a Lipschitz map, where $U\subset\R^k$ is open. 
The metric directional derivative of $\varphi$ in direction $v\in\R^k$ is defined by
\begin{equation*}
 \md\varphi_z(v):= \lim_{r\searrow 0}\frac{d(\varphi(z+rv),\varphi(z))}{r}
\end{equation*}
if this limit exists. 
It was proved independently by Kirchheim \cite{Kirchheim} and Korevaar-Schoen \cite{Korevaar-Schoen} that for almost every $x\in U$ the metric
derivative $\md\varphi_x(v)$ exists for all $v\in\R^k$ and defines a seminorm on $\R^k$. The following theorem is a consequence of this metric differentiability property.

\bt\label{theorem:kirchheim-bilip-finite-set}
 Let $(X,d)$ be a metric space and $\varphi:K\to X$ Lipschitz with $K\subset\R^k$ Borel measurable and such that $\hm^k(\varphi(K))>0$. Then 
 there exists a norm $\|\cdot\|$ on $\R^k$ with the following property: For every $\varepsilon>0$ and for every finite set $\Lambda\subset\R^k$ there exist
 $r>0$ and a map $\psi: \Lambda\to X$ such that $\psi: (\Lambda,r\|\cdot\|)\to X$ is $(1+\varepsilon)$-biLipschitz.
\et

For the proof see Lemma 4 and Theorem 7 (area formula) of \cite{Kirchheim}. The norm $\|\cdot\|$ in \thmref{theorem:kirchheim-bilip-finite-set} is in fact given by $\|\cdot\|=\md\varphi_x$
for a suitable $x\in K$.

The Jacobian of a seminorm $s$ on $\R^k$ is defined by
\begin{equation*}
 \jac_k(s):= \frac{\omega_k}{\lm^k(\{v\in\R^k: s(v)\leq 1\})},
\end{equation*}
where $\omega_k$ is the volume of the unit ball in $\E^k$ and $\lm^k$ is the Lebesgue measure. If $\varphi: D^2\to X$ is Lipschitz then its parametrized Hausdorff area is
\begin{equation*}
 \area(\varphi):= \int_{D^2}\jac_2(\md\varphi_x)d\lm^2(x).
\end{equation*}
If $\varphi$ is one-to-one on a subset of $D^2$ of full measure then, by Theorem 7 of \cite{Kirchheim}, we have $\area(\varphi)=\hm^2(\varphi(D^2))$.

\subsection{Integral currents in metric spaces}\label{section:currents}
The theory of integral currents in metric spaces was developed by Ambrosio and Kirchheim in \cite{Ambr-Kirch-curr} and provides a suitable 
notion of surfaces and area/volume in the setting of metric spaces. In the following we adopt the notation of \cite{Ambr-Kirch-curr} and refer to it for precise definitions. The
definitions which will be needed throughout this article can also be found in Section 2.2 of \cite{Wenger-GAFA}.

For $k\geq 0$ the space of $k$-dimensional integral currents in a complete metric space $X$ is denoted by $\intcurr_k(X)$, the mass measure of an 
element $T\in\intcurr_k(X)$ by $\|T\|$ and its mass by $\mass{T}:=\|T\|(X)$. If $k\geq 1$, the boundary of $T$ is denoted by $\bdry T$ and is an element of $\intcurr_{k-1}(X)$.
It will be shown in \lemref {lemma:almost-optimal-decomposition-one-cycles} that $1$-dimensional integral currents are essentially induced by Lipschitz curves. As regards
$2$-dimensional integral currents, an element $T\in\intcurr_2(X)$ can be thought of as a $2$-dimensional oriented surface (with arbitrary genus and possibly with integer multiplicity) 
which is locally parametrized by biLipschitz maps from $\R^2$ and whose boundary consists of a union of Lipschitz curves of finite total length. Moreover, $\|T\|$ is a particular Finsler 
area on the surface taken with multiplicity, namely Gromov's $2$-dimensional mass$*$ area defined in \cite{Gromov-filling}. Of course, in the setting of Riemannian manifolds this is 
simply the Riemannian area. Singular Lipschitz discs and singular Lipschitz chains in $X$ (in the sense of Gromov \cite{Gromov-filling}) induce in a natural way $2$-dimensional integral currents.

The following definitions and constructions are frequently used throughout this text. 
Every Borel subset $A\subset\R^k$ with finite measure and finite perimeter induces an element of $\intcurr_k(\R^k)$ by
\begin{equation*}
 \Lbrack\chi_A\Rbrack(f,\pi_1,\dots,\pi_k):= \int_Af\det\left(\frac{\partial\pi_i}{\partial x_j}\right)\,d\lm^k.
\end{equation*}
Recall that (integral) currents of dimension $k$ in $X$ are in particular functionals on the space of $(k+1)$-tuples $(f,\pi_1,\dots,\pi_k)$ of Lipschitz functions on $X$ with $f$ bounded.
Given $T\in\intcurr_k(X)$ and a Lipschitz map $\varphi: X\to Y$, where $Y$ is another complete metric space, the pushforward of $T$ by $\varphi$ is defined by
\begin{equation*}
 \varphi_{\#}T(g,\tau_1,\dots,\tau_k):= T(g\circ\varphi,\tau_1\circ\varphi,\dots,\tau_k\circ\varphi)
\end{equation*}
and is an element of $\intcurr_k(Y)$. It can be shown that $\mass{\varphi_{\#}T}\leq \lip(\varphi)^k\mass{T}$, where $\lip(\varphi)$ is the Lipschitz constant of $\varphi$.
The boundary of $T\in\intcurr_k(X)$ is
\begin{equation*}
 \bdry T(f,\pi_1,\dots\pi_{k-1}):= T(1,f,\pi_1,\dots,\pi_{k-1})
\end{equation*}
and defines an element of $\intcurr_{k-1}(X)$. It follows directly from the definitions that $\bdry(\varphi_{\#}T) = \varphi_{\#}(\bdry T)$.
A Lipschitz curve $\gamma:[a,b]\to X$ gives rise to the element $\gamma_{\#}\Lbrack\chi_{[a,b]}\Rbrack\in\intcurr_1(X)$ where $\chi_{[a,b]}$ denotes the characteristic function.
If $\varphi$ is one-to-one then 
$\mass{\gamma_{\#}\Lbrack\chi_{[a,b]}\Rbrack}=\length(\gamma)$. 
If $\gamma$ is a Lipschitz loop and $S\in\intcurr_2(X)$ satisfies $\bdry S=\gamma_{\#}\Lbrack\chi_{[a,b]}\Rbrack$ then $\gamma$ is said to bound $S$.
A Lipschitz map $\varphi:D^2\to X$ gives rise to the $2$-dimensional integral current $S:= \varphi_{\#}\Lbrack\chi_{D^2}\Rbrack$. It should be noted that
in general $\area(\varphi)\not=\mass{S}$, since $\mathbf{M}$ corresponds to the mass$*$ area rather than the Hausdorff area, see also \secref{section:area-functionals}.
A singular Lipschitz chain $c=\sum m_i\varphi_i$ gives rise to the integral current $\sum m_i\varphi_{i\#}\Lbrack\chi_{\Delta}\Rbrack$.

The following lemma will be needed in the proof of \propref{proposition:tree-bilip}. It shows that integral $1$-currents without boundary are essentially countable unions of Lipschitz
loops.
\bl\label{lemma:almost-optimal-decomposition-one-cycles}
 Let $X$ be a complete metric length space, $T\in\intcurr_1(X)$ with $\bdry T=0$ and $\varepsilon>0$. There then exist at most countably many Lipschitz 
 loops $\gamma_i:[0,a_i]\to X$ with the properties that
 \begin{equation}\label{equation:current-is-lipschitz-curve}
  T = \sum_i\gamma_{i\,\#}\Lbrack\chi_{[0,a_i]}\Rbrack
 \end{equation}
 and
 \begin{equation*}
  \sum_i\length(\gamma_i)\leq (1+\varepsilon)\mass{T}.
 \end{equation*}
\el

Note that \eqref{equation:current-is-lipschitz-curve} implies that 
\begin{equation*}
 \mass{T}\leq \sum\length(\gamma_i)
\end{equation*}
by the sub-additivity of $\mathbf M$.
In Euclidean space it can be shown (see \cite[4.2.25]{Federer}) that \lemref{lemma:almost-optimal-decomposition-one-cycles} holds with $\varepsilon=0$,
and it is conceivable that the same should be true in all complete length spaces.

\begin{proof}
 It clearly suffices to prove that for every $\varepsilon>0$ there exist finitely many Lipschitz loops $\gamma_i:[0,a_i]\to X$ such that
 \begin{equation*}
  \mass{T-\sum \gamma_{i\#}\Lbrack\chi_{[0,a_i]}\Rbrack}\leq\varepsilon\mass{T}
 \end{equation*}
 and
 \begin{equation*}
  \sum\length(\gamma_i)\leq (1+\varepsilon)\mass{T}.
\end{equation*}
 In order to find such a decomposition let $\varepsilon'>0$ be small enough, to be determined later. Using Lemma 4 and Theorem 7 of \cite{Kirchheim} one easily shows that
 there exist finitely many $(1+\varepsilon')$-biLipschitz maps $\varphi_i: K_i\to X$, $i=1,\dots, n$, where 
 $K_i\subset\R$ are compact and such that $\varphi_i(K_i)\cap\varphi_j(K_j)=\emptyset$ if $i\not=j$, and 
 \begin{equation*}
  \|T\|\left(X\ohne \cup\varphi_i(K_i)\right)\leq\varepsilon'\mass{T},
 \end{equation*}
 see also \cite[Lemma 4.1]{Ambr-Kirch-curr}.
 By McShane's extension theorem there exists a $(1+\varepsilon')$-Lipschitz extension $\overline{\eta}_i: X\to\R$ of $\varphi_i^{-1}$ for each $i=1,\dots,n$. 
 Set $\Omega:= \cup\varphi_i(K_i)$ and let $\{z_1,\dots, z_m\}\subset\Omega$ be a finite and $\delta$-dense set for $\Omega$, where $\delta>0$ is such that
 \begin{equation}\label{equation:equation:condition-K-distance}
  \dist(\varphi_i(K_i),\varphi_j(K_j))\geq\frac{2(1+\varepsilon')\delta}{\varepsilon'}\quad\text{whenever $i\not=j$.}
 \end{equation}
 We set $N:= m+n$ and define a map $\Psi: X\to \ell^\infty_N$ by
 \begin{equation*}
  \Psi(x):= \left(\overline{\eta}_1(x),\dots,\overline{\eta}_n(x), d(x,z_1),\dots,d(x,z_m)\right).
 \end{equation*}
 Note that $\Psi$ is $(1+\varepsilon')$-Lipschitz and $(1+\varepsilon')$-biLipschitz on $\Omega$. Indeed, it is clear that the latter statement holds when
 restricted to each $\varphi_i(K_i)$. 
 Moreover, for $x\in \varphi_i(K_i)$ and $x'\in \varphi_j(K_j)$ with $i\not=j$ there exists a $z\in\varphi_i(K_i)$ with $d(x,z)\leq\delta$ and hence
 \begin{equation*}
  d(x,x')\leq d(x,z)+d(z,x')\leq\|\Psi(x')-\Psi(x)\|_\infty+2\delta
 \end{equation*}
 from which the biLipschitz property follows together with \eqref{equation:equation:condition-K-distance}.
 By \cite[4.2.25]{Federer} there exist Lipschitz curves $\varrho_j:[0,a_j]\to\ell^\infty_N$ which are parametrized by arc-length, one-to-one on $(0,a_j)$, 
 with $\varrho_j(0)=\varrho_j(a_j)$ and satisfy $\Psi_\#T= \sum_{j=1}^\infty\varrho_{j\#}\Lbrack\chi_{[0,a_j]}\Rbrack$ and
 \begin{equation*}
  \mass{\Psi_\#T}= \sum_{j=1}^\infty\mass{\varrho_{j\#}\Lbrack\chi_{[0,a_j]}\Rbrack}=\sum_{j=1}^\infty a_j.
 \end{equation*}
 Choose $M\in\N$ sufficiently large such that $R:= \sum_{j=M+1}^\infty \varrho_{j\#}\Lbrack\chi_{[0,a_j]}\Rbrack$ satisfies
 \begin{equation*}
  \mass{R}\leq\varepsilon'\mass{T}.
 \end{equation*}
 Since $X$ is a length space there exists a $(1+2\varepsilon')$-Lipschitz extension $\gamma_j:[0,a_j]\to X$ of 
 $(\Psi\on{\Omega})^{-1}\circ(\varrho_j\on{\varrho_j^{-1}(\Psi(\Omega))})$ with $\gamma_j(a_j)=\gamma_j(0)$ for each $j=1,\dots, M$. 
 We now have
  \begin{equation*}
  \sum_{j=1}^M \varrho_{j\#}\Lbrack\chi_{\varrho_j^{-1}(\Psi(\Omega)^c)}\Rbrack = \left[\Psi_\#(T\rstr\Omega^c)- R\right]\rstr\Psi(\Omega)^c
 \end{equation*}
 from which it easily follows that
 \begin{equation*}
  \begin{split}
   T-\sum_{j=1}^M&\gamma_{j\#}\Lbrack\chi_{[0,a_j]}\Rbrack \\ &= \left(\Psi\on{\Omega}\right)^{-1}_\#\left[(R-\Psi_\#(T\rstr\Omega^c))\rstr\Psi(\Omega)\right]
   - \sum_{j=1}^M\gamma_{j\#}\Lbrack\chi_{\varrho_j^{-1}(\Psi(\Omega)^c)}\Rbrack
   + T\rstr\Omega^c
  \end{split}
 \end{equation*}
 and
 \begin{equation*}
  \sum_{j=1}^M\hm^1(\varrho_j^{-1}(\Psi(\Omega)^c))=\sum_{j=1}^M\mass{\varrho_{j\#}\Lbrack\chi_{\varrho_j^{-1}(\Psi(\Omega)^c)}\Rbrack}
   \leq \varepsilon'(2+\varepsilon')\mass{T}.
 \end{equation*}
 This leads to
 \begin{equation*}
  \mass{T-\sum_{j=1}^M\gamma_{j\#}\Lbrack\chi_{[0,a_j]}\Rbrack}\leq [5+8\varepsilon'+3\varepsilon'^2]\varepsilon'\mass{T}.
 \end{equation*}
 Finally, we estimate
 \begin{equation*}
  \sum_{j=1}^M\length(\gamma_j)\leq (1+2\varepsilon')\sum_{j=1}^Ma_j\leq(1+2\varepsilon')\mass{\Psi_\#T}\leq (1+2\varepsilon')(1+\varepsilon')\mass{T}.
 \end{equation*}
 This proves the claim at the beginning of the proof given that $\varepsilon'>0$ was chosen small enough.
\end{proof}

\subsection{Area functionals and the isoperimetrix}\label{section:area-functionals}
In normed spaces various definitions of area and volume have been studied, see e.g.\ the survey \cite{Alvarez-Thompson}. 
These definitions can be used to define area and volume functionals also for integral currents, as is explained
in Section 13 of \cite{Ambr-Kirch-curr}. It turns out that our results hold for various definitions of area. The facts below will be needed in the proofs of
\thmref{theorem:isop-special-case} and of the results in \secref{section:isop-gromov}.

Fix a definition of area $\mu$, see \cite{Alvarez-Thompson} for this terminology. Thus, $\mu$ assigns to every $2$-dimensional normed space $V$ a Haar measure $\mu_V$ on $V$
(in particular, $\mu_V=c_V\hm^2$, where $c_V$ varies continuously with $V$ and $c_V=1$ if $V=\E^2$). 
We denote by $\mu^b$ the Hausdorff, by $\mu^{ht}$ the Holmes-Thompson, and by $\mu^{m*}$ the Gromov mass$*$ definition of area.
Let $\mathbb{I}_V$ be an isoperimetric subset of $V$, that is a compact convex subset of $V$ with non-empty interior which has maximal $\mu_V$-area among all subsets with the 
same perimeter.
It can be shown that if $\mu$ is one of the three area definitions above then
\begin{equation*}
 \mu(\mathbb{I}_V)\geq \frac{1}{4\pi}\length(\bdry\mathbb{I}_V)^2,
\end{equation*}
see e.g.~\cite[p.~33]{Alvarez-Thompson}. In fact, if $\mu=\mu^{ht}$ then we always have equality. On the other hand, if $\mu=\mu^{m*}$ then 
equality holds if and only if $V$ is the Euclidean plane.
It can furthermore be shown (see Theorem 3.13 of \cite{Alvarez-Thompson}) that $\mu^{ht}_V\leq\mu^b_V$ and $\mu^{ht}_V\leq\mu^{m*}_V$ 
for all ($2$-dimensional) normed spaces $V$. 

For the following we refer to Section 13 of \cite{Ambr-Kirch-curr}. Denote by $\mathcal{F}_\mu$ the area functional for integral currents induced by $\mu$. For example, $\mathcal{F}_{\mu^{m*}}$ is simply the mass $\mathbf M$. Furthermore,
if $\varphi:D^2\to X$ is Lipschitz then $\mathcal{F}_{\mu^b}(\varphi_{\#}\Lbrack\chi_{D^2}\Rbrack)=\area(\varphi(D^2))$. It should be noted that the area functionals associated to 
$\mu^b$, $\mu^{ht}$ and $\mu^{m*}$ all agree up to a universal constant (by John's theorem). In the proof of \thmref{theorem:isop-general-case} we will need
the following (well-known and easy to prove) 
semi-ellipticity property of $\mathbf{M}$: Let $W$ be a normed space, $V\subset W$ a $2$-dimensional affine subspace and $C\subset V$
a compact convex set. Then $\mu^{m*}(C)\leq \mass{S}$ for every $S\in\intcurr_2(W)$ whose boundary is induced by a Lip\-schitz loop which parametrizes $\bdry C$.
See \cite[Theorem 4.28]{Alvarez-Thompson} for a much stronger statement. In the recent major advance \cite{Burago-Ivanov} it has been shown that
$\mu^{ht}(C)\leq \mathcal{F}_{\mu^{ht}}(\Sigma)$ for all singular Lipschitz discs with boundary $\bdry C$. It is not known whether the same holds if $\Sigma$ is replaced by
a singular Lipschitz chain.
Furthermore, it is a long-standing open question going back to Busemann whether an analogous statement holds for $\mu^b$.

We end this section with the following simple but crucial fact.

\bl \label{lemma:Linfty-fillings}
 Let $X$ and $Y$ be metric spaces and suppose $X$ isometrically embeds in $Y$. Let $\mu$ be a definition of volume, $k\geq 1$, and
 $T\in\intcurr_k(X)$ with $\bdry T=0$. Then for every $S\in\intcurr_{k+1}(Y)$ with $\bdry S = T$ there exists $S'\in\intcurr_{k+1}(L^\infty(X))$ with $\bdry S'=T$ and such that
 \begin{equation*}
  \mathcal{F}_{\mu}(S')\leq\mathcal{F}_{\mu}(S).
 \end{equation*}
\el
This follows indeed directly from the Lipschitz extension property of $L^\infty(X)$ and the fact that $\mathcal{F}_{\mu}(\varphi_{\#}S)\leq \mathcal{F}_{\mu}(S)$ if 
$\varphi$ is $1$-Lipschitz.

\section{Isoperimetric inequalities of thickenings}\label{section:coarse-isoperimetric}

Let $X$ and $X'$ be metric spaces.  $X'$ is called a thickening of $X$ if there exists an isometric embedding $\varphi: X\to X'$ such that $\varphi(X)$ is in 
finite Hausdorff distance of $X'$. A complete metric space $Y$ is said to admit a quadratic isoperimetric inequality for curves if there exists $C>0$ 
such that every Lipschitz loop $\gamma$ in $Y$ bounds an $S\in\intcurr_2(Y)$ with
\begin{equation*}
 \mass{S}\leq C\length(\gamma)^2.
\end{equation*}
In contrast, $Y$ is said to admit a quadratic isoperimetric inequality for sufficiently long curves if the above holds for all $\gamma$ with $\length(\gamma)\geq s_0$ for some
$s_0\geq 0$.
%

\bd\label{definition:admissible}
 A metric space $X$ is called admissible if there exists a complete metric space $X_\delta$ which is a thickening of $X$ and which admits a quadratic isoperimetric 
 inequality for curves.
\ed

In the main result of this section, \propref{proposition:admissible-thickening}, we will show that metric length spaces with a coarse quadratic isoperimetric inequality or a quadratic
isoperimetric inequality for sufficiently long curves are admissible. As mentioned in the introduction the latter will be needed in the proof of \thmref{theorem:isop-special-case}.

The following notion of coarse homological fillings of Lipschitz loops in length spaces generalizes that of coarse fillings by discs given in \cite[III.H.2]{Bridson-Haefliger}.
Let $X$ be a length space, $\delta, a>0$ and $\gamma:[0,a]\to X$ a Lipschitz loop. A $\delta$-coarse
homological filling of $\gamma$ is a triple $(K, c, \mu)$ with the following properties:
\begin{enumerate}
 \item $K$ is a $2$-dimensional simplicial complex such that every attaching map of a $2$-cell is a $3$-gon;
 \item $c$ is a simplicial integral $2$-chain in $K$, that is a function on the $2$-cells with values in $\Z$;
 \item $\mu: K\to X$ is a (possibly discontinuous) map such that
  \begin{enumerate}
   \item $\diam(\mu(e))\leq \delta$ for each $2$-cell $e\subset K$;
   \item there exists a combinatorial map $\varrho:[0,a]\to K^{(1)}$ such that $\gamma = \mu\circ\varrho$ and such that the $1$-cycle induced by $\varrho$ is $\bdry c$.
    Here, $[0,a]$ is endowed with a combinatorial structure of the form $$[0,a]= [s_0, s_1]\cup\dots\cup[s_{n-1}, s_n]$$ for some $0=s_0<\dots<s_n=a$.
  \end{enumerate}
\end{enumerate}
For the definition of simplicial $2$-complexes we refer to \cite[I.8A.4]{Bridson-Haefliger}. We mention here that every simplex of dimension $1$ or $2$ in $K$ inherits an 
orientation coming from its attaching map. It is clear that every $\delta$-coarse filling as defined in \cite[III.H.2]{Bridson-Haefliger} induces a $\delta$-coarse homological 
filling.
In \cite{Gersten-subgroups} Gersten introduced a homological notion of fillings in the context of simplicial complexes and groups. His approach
uses surface diagrams the foundation of which were laid in the book \cite{Olshanskii-relation}.
In the sequel a $\delta$-coarse homological filling will simply be called a $\delta$-filling. If $K$ is homeomorphic to a disc then the filling will be called a $\delta$-coarse disc filling.
The $\delta$-area of the triple $(K,c,\mu)$ is by definition 
\begin{equation*}
\area_\delta(K,c,\mu):= \sum_{e\text{ $2$-cell in $K$}}  |c(e)|
\end{equation*}
and the $\delta$-filling area of $\gamma$ is given by
\begin{equation*}
 \fillarea_\delta(\gamma):= \inf\left\{\area_\delta(K,c,\mu): (K,c,\mu)\text{ $\delta$-filling of $\gamma$}\right\}.
\end{equation*}
A function $f$ for which
\begin{equation*}
 \fillarea_\delta(\gamma)\leq f(\length(\gamma))
\end{equation*}
for every Lipschitz loop $\gamma: [0,a]\to X$ is called a $\delta$-coarse homological isoperimetric bound for $X$.

It is not difficult to establish the following quasi-isometry invariance of coarse homological isoperimetric bounds. The proposition will not be used anywhere in the text and 
its proof is omitted.
\bp
 Let $X$ and $Y$ be quasi-isometric length spaces and suppose $X$ has a $\delta$-coarse homological isoperimetric bound $f$ for some $\delta>0$.  
 Then there exists $\delta'>0$ such that $Y$ has a $\delta'$-coarse homological isoperimetric bound $g$ satisfying $g\preceq f$.
\ep
Here, $g\preceq f$ means that there exists $K>0$ such that $g(s)\leq Kf(Ks+K)+ Ks + K$ for all $s\geq 0$.
The proof of an analogous statement for $\delta$-coarse disc fillings can be found in \cite[III.H.2.2]{Bridson-Haefliger}.

The following is the main result of this section.

\bp\label{proposition:admissible-thickening}
 Let $X$ be a length space. If $X$ admits a coarse homological quadratic isoperimetric inequality then $X$ is admissible in the sense of \defref{definition:admissible}.
 Similarly, if $X$ admits a quadratic isoperimetric inequality for sufficiently long curves then $X$ is admissible.
\ep

As an immediate consequence we obtain that Cayley graphs of finitely presented groups with at most quadratic Dehn function as well as metric spaces satisfying the hypotheses of 
\thmref{theorem:isop-special-case} are admissible. Of course, geodesic Gromov hyperbolic spaces are all admissible, they have geodesic thickenings even with a linear isoperimetric 
inequality for curves.

\begin{proof}
 Denote by $\overline{X}$ the completion of $X$. Let $Z\subset\overline{X}$ be a maximally $\delta$-separated subset. For $z\in Z$ denote by 
 $X_z$ the injective envelope of $B_z:=(B(z,8\delta), d_{B(z,8\delta)})$, where $d_{B(z,8\delta)}$ is the length metric on the ball $B(z,8\delta)\subset\overline{X}$.
 Denote the metric on $X_z$ by $d_z$. First of all, it is clear that $B_z$ is closed in $X_z$. Furthermore, one can easily show that 
 $\diam X_z\leq 64\delta$.
 Now set
 \begin{equation*}
  X_\delta:=\Big(\bigsqcup_{z\in Z}X_z\Big)_{\;\;\big/\;{ \sim}}
 \end{equation*}
 where $x\sim x'$ if and only if $x\in B_z\subset X_z$ and $x'\in B_{z'}\subset X_{z'}$ for some $z,z'\in Z$ and $x=x'$.
 Define a metric on $X_\delta$ as follows. For $x\in X_z$ and $x'\in X_{z'}$ set 
  \begin{equation*}
  \varrho_{zz'}(x,x'):=
    \inf\left\{d_z(x,y)+d_X(y,y')+d_{z'}(y', x'): y\in B_z, y'\in B_{z'}\right\}
 \end{equation*}
 and
 \begin{equation*}
  d_\delta([x],[x']):= \left\{
   \begin{array}{c@{\quad:\quad}l}
    \varrho_{zz'}(x,x') & z\not=z'\\
    \min\{\varrho_{zz}(x,x'), d_z(x,x')\} & z=z'.
    \end{array}\right.
 \end{equation*}
 It is straight-forward to check that $d_\delta$ defines a complete metric on $X_\delta$ and that $(X_\delta, d_\delta)$ is a length space which isometrically contains 
 $\overline{X}$ as a closed subset.
 Moreover, $X_\delta$ is a thickening of $X$.

 Let now $a>0$ and set $b:=\min\{a, 1\}$. Let $\gamma:[0,a]\to X_\delta$ be a Lipschitz loop and assume $\gamma([0,a])$ is not entirely contained in a single $X_z$ and that 
 $\gamma(0)\in\overline{X}$.
 We construct a Lipschitz map $\varphi: Q\to X_\delta$, where $Q:=[0,a]\times [0,b]$, with the property that 
 \begin{equation*}
  \varphi(s,t) = \left\{\begin{array}{l@{\quad:\quad}l}
    \gamma(s) & t=0\\
    \gamma(0) & s\in \{0,a\}\\
    \psi(s) & t=b
   \end{array}\right.
 \end{equation*}
 for some $2\lip(\gamma)$-Lipschitz loop $\psi:[0,a]\to X$ of length at most $2\length(\gamma)$.
 For this set $U:= \gamma^{-1}(X_\delta\ohne\overline{X})$. We may assume without loss of generality that $U\not=\emptyset$.
 For each $u\in U$ let $\tau_u$ and $\nu_u$ be the smallest and largest value, respectively, such that 
 $U_u:= (\tau_u, \nu_u)\subset U$
 and $u\in U_u$. Note that for $u,u'\in U$ we either have $U_u=U_{u'}$ or $U_u\cap U_{u'}=\emptyset$. Therefore, there exist countably many $u_j\in U$ such that
 $U=\sqcup U_{u_j}$. For each $j\in\N$ let $z_j\in Z$ be such that $\gamma(u_j)\in X_{z_j}$. By construction there exist Lipschitz curves
 $\tilde{\gamma}_j: [\tau_{u_j}, \nu_{u_j}]\to B_{z_j}$ parametrized proportionally to arc-length  and such that 
 $\tilde{\gamma}_j(\tau_{u_j})=\gamma(\tau_{u_j})$
 and  $\tilde{\gamma}_j(\nu_{u_j})=\gamma(\nu_{u_j})$ and $\length(\tilde{\gamma}_j)\leq 2\length(\gamma\on{U_j})$. 
 Define $\psi: [0,a]\to\overline{X}$ by setting 
 $\psi(s):= \tilde{\gamma}_j(s)$ if $s\in U_{u_j}$ and $\psi(s):= \gamma(s)$ otherwise and note that $\psi$ is a Lipschitz 
 loop with Lipschitz constant at most $2\lip(\gamma)$ and of length at most $2\length(\gamma)$.
 Set $\varphi(s,0):= \gamma(s)$ and $\varphi(s,b):=\psi(s)$ for all $s\in[0,a]$.
 Furthermore, define $\varphi(\tau_{u_j},t):= \gamma(\tau_{u_j})$ and $\varphi(\nu_{u_j},t):= \gamma(\nu_{u_j})$
 for all $j$ and all $t\in[0,b]$. Now, it is clear from the injectivity of every $X_z$ that $\varphi$ can be extended to a Lipschitz map
 $\varphi: Q\to X_\delta$ with Lipschitz constant at most $C_1\lip(\gamma)$ for some $C_1$ depending only on $\delta$.  
 Using again the injectivity of the $X_z$ and the fact that $Z$ is $2\delta$-dense in $\overline{X}$ it can easily be seen that $\psi$ can in fact be constructed so that its image 
 is in $X$ instead in $\overline{X}$.

In the following, $C_2, C_3$, and $C_4$ will denote constants only depending on $\delta$.
Let $\gamma:[0,a]\to X_\delta$ be a Lipschitz loop and let $T$ be the integral current induced by $\gamma$, that is $T:=\gamma_{\#}\Lbrack\chi_{[0,a]}\Rbrack$. 
We may assume that $\gamma$ is parametrized by arc-length, thus $a=\length(\gamma)$. 
 If $\gamma([0,a])$ is contained in a single $X_z$ then there exists an $C_2$-Lipschitz map $\varphi:Q\to X_z$ with $\varphi(s,0)=\gamma(s)$ and 
 with $\varphi(s,t)=\gamma(0)$ if $s\in\{0,a\}$ or $t=b$. Consequently, $S:=\varphi_\#\Lbrack\chi_Q\Rbrack$ satisfies $S\in\intcurr_2(X_\delta)$, 
 $\bdry S = T$ and $\mass{S}\leq C_2^2ab$.
 On the other hand, if $\gamma([0,a])$ is not contained in a single $X_z$ then we may assume that $\gamma(0)\in \overline{X}$ after a possible change of parametrization.
 By the above there then exists a $C_1$-Lipschitz map $\varphi:Q\to X_\delta$ with $\varphi(s,0)=\gamma(s)$ and
 $\varphi(s,b)=\psi(s)$ for all $s\in[0,a]$, where $\psi: [0,a]\to X$ is a $2$-Lipschitz loop of length at most $2a$, and
 $\varphi(s,t)=\gamma(0)$ if $s\in\{0,a\}$.
 Let $(K,c,\mu)$ be a $2\delta$-coarse filling of $\psi$ in $X$ with $\area_\delta(K, c, \mu)\leq 4Ca^2$. Let furthermore $\varrho:[0,a]\to K^{(1)}$ be as in the 
 definition of the coarse filling. Construct a
 map $\overline{\mu}: K\to X_\delta$ as follows. Set $\overline{\mu}(z):= \mu(z)$ whenever $z\in K^{(0)}$ and, in a first step, extend $\overline{\mu}$ 
 to $K^{(1)}$ in such a way that $\overline{\mu}\on{e}$ is a Lipschitz curve in $X$ parametrized proportional to arc-length 
 joining its endpoints and of length at most $\frac{3}{2}\diam(\mu(\bdry e))\leq 3\delta$, for each $1$-cell $e\subset K^{(1)}$.
 Hereby, each $e$ is to be induced with the Euclidean metric.
 Since for every closed $2$-cell $e\subset K$ we have $\diam(\overline{\mu}(\bdry e))\leq 6\delta$ and since $Z$ is $2\delta$-dense in $X$ we obtain that 
 $\overline{\mu}(\bdry e)$ is contained in $B(z,8\delta)$ for some $z\in Z$. Furthermore, $\overline{\mu}\on{\bdry e}: \bdry e\to B_z$ is $C_3\delta$-Lipschitz 
 and hence can be extended to a $C_3\delta$-Lipschitz map $\overline{\mu}\on{e}:e\to X_z$. This yields the desired map 
 $\overline{\mu}$. 
 We can moreover construct a $C_4$-Lipschitz map $\overline{\psi}: Q\to X_\delta$ with the property that
 \begin{equation*}
  \overline{\psi}(s,t) = \left\{\begin{array}{l@{\quad:\quad}l}
    \psi(s) & t=0\\
    \overline{\mu}\circ\overline{\varrho}(s) & t=b\\
    \psi(s_i) & s=s_i.
   \end{array}\right.
 \end{equation*}
 Here, $\overline{\varrho}:[0,a]\to K^{(1)}$ is a reparametrization of $\varrho$ such that $\overline{\varrho}\on{[s_i,s_{i+1}]}$ is a constant-speed 
 parametrization of the $1$-cell $\varrho([s_i,s_{i+1}])$ for all $i$.
 As for this construction it is enough to note that we have $\psi(s_i)=\overline{\mu}\circ\overline{\varrho}(s_i)$ for all $i=0,\dots,n$ and that 
 \begin{equation*}
  \diam\left(\psi([s_i, s_{i+1}])\cup\overline{\mu}(\overline{\varrho}([s_i,s_{i+1}]))\right)\leq 5\delta.
 \end{equation*}
 The existence of $\overline{\psi}$ now follows from the same arguments as above.

 Finally, we can define a suitable filling of $T$ by setting 
 \begin{equation*}
  S:= \varphi_\#\Lbrack\chi_Q\Rbrack + \overline{\psi}_\#\Lbrack\chi_Q\Rbrack  + \overline{\mu}_\#c.
 \end{equation*}
 Indeed, we have $S\in\intcurr_2(X_\delta)$ and 
 \begin{equation*}
   \bdry S = \varphi_\#(\bdry\Lbrack\chi_Q\Rbrack) + \psi_\#\Lbrack\chi_{[0,a]}\Rbrack
     - (\overline{\mu}\circ\overline{\varrho})_\#\Lbrack\chi_{[0,a]}\Rbrack
     + \overline{\mu}_\#(\bdry c) = T
 \end{equation*}
 as well as
 \begin{equation*}
   \mass{S} \leq \lip(\varphi)^2ab + \lip(\overline{\psi})^2ab + 4C'C_3^2\delta^2 a^2
 \end{equation*}
 for a constant $C'$ only depending on $C$.
 This completes the proof of the first statement. The second statement uses the same constructions above.
\end{proof}

\section{Asymptotic subsets and Gromov hyperbolicity}\label{section:asymptotic-subsets}
A metric space $(Z,d_Z)$ is said to be an asymptotic subset of another metric space $(X,d_X)$ if
there exist a sequence of subsets $Z_n\subset X$ and $r_n\nearrow\infty$ such that $(Z_n,r_n^{-1}d_X)$ converges in the Gromov-Hausdorff sense to $(Z,d_Z)$. 

The proposition below, the main result of this section, plays a crucial role in the proof of our main theorem.

\bp\label{proposition:tree-bilip}
 Let $X$ be an admissible geodesic metric space and suppose that $\hm^2(\varphi(K))=0$ whenever $\varphi: K\to (Z,d_Z)$ is a Lipschitz map
 with $K\subset\R^2$ compact and $(Z,d_Z)$ an asymptotic subset of $X$. Then $X$ is Gromov hyperbolic.
\ep

It should be noted that asymptotic subsets can be replaced by asymptotic cones in the above statement. 
For the proof we will need the following construction.
Given a geodesic metric space and $\lambda>0$ define a function by
\begin{equation}\label{equation:non-degeneracy-integral}
 H_\lambda(r):= \sup_{f,\pi,\gamma}\int_0^1(f\circ\gamma)(s)(\pi\circ\gamma)'(s)ds,
\end{equation}
where the supremum is taken over all $\lambda r$-Lipschitz curves $\gamma:[0,1]\to X$ with $\gamma(1)=\gamma(0)$ and $\lambda r^{-1}$-Lipschitz functions $f,\pi: X\to\R$.
 Note that the integral in \eqref{equation:non-degeneracy-integral} corresponds exactly to  $T(f,\pi)$, where $T$ is the integral current given by $T:=\gamma_\#\Lbrack\chi_{[0,1]}\Rbrack$, and $T(f,\pi)$ remains 
 unchanged when $f,\pi$ are replaced by $f+c_1$ and $\pi+c_2$ for constants $c_1,c_2$.
 The definitions of $\bdry S$ and $\|S\|$, see (2.2) of \cite{Ambr-Kirch-curr} for the latter, furthermore yield
 \begin{equation*}
 |T(f,\pi)|=|S(1,f,\pi)|\leq \lip(f)\lip(\pi)\mass{S}
 \end{equation*}
 whenever $S\in\intcurr_2(L^\infty(X))$ satisfies $\bdry S=T$ and thus
 \begin{equation}\label{equation:H-basic-estimate}
  \left|\int_0^1(f\circ\gamma)(s)(\pi\circ\gamma)'(s)ds\right|\leq \frac{\lambda^2}{r^2}\fillarea_{L^\infty(X)}(\gamma),
 \end{equation}
 where $\fillarea_{L^\infty(X)}(\gamma)$ is the least mass of an $S'\in\intcurr_2(L^\infty(X))$ with $\bdry S'=\gamma_\#\Lbrack\chi_{[0,1]}\Rbrack$.
 In particular, we obtain that $H_\lambda(r)\leq \lambda^4$ by the cone inequality \cite[Proposition 10.2]{Ambr-Kirch-curr}.
The function $H_\lambda$ in some sense measures how `collapsed' closed curves in $X$ are. 

\bl\label{lemma:degenerate-loops}
 A geodesic metric space $X$ is Gromov hyperbolic if and only if 
 \begin{equation}\label{equation:lim-H-0}
 \lim_{r\to\infty}H_\lambda(r)=0\quad\text{ for every $\lambda>0$.}
 \end{equation}
\el

\begin{proof}
 We first prove by contradiction that \eqref{equation:lim-H-0} implies Gromov hyperbolicity. 
 Assume therefore that $X$ is not Gromov hyperbolic. By the main theorem of \cite{Bonk-detour} there exists $C\in(0,\infty)$, a sequence $r_n\nearrow\infty$ and curves
 $\alpha_n:[0,1]\to X$ of length bounded above by $Cr_n$ with the following property: For every $n\in\N$ there exists a geodesic segment $\beta_n:[0,1]\to X$ from $\alpha_n(0)$ to 
 $\alpha_n(1)$ such that
 \begin{equation*}
  \image(\alpha_n)\cap B(z_n,r_n) = \emptyset
 \end{equation*}
 for some $z_n\in\image(\beta_n)$.
 Let $\gamma_n:[0,1]\to X$ be the concatenation of $\beta_n$ and $\alpha_n$, parametrized proportional to arc-length. Define
 \begin{equation*}
  f_n(x):= \max\{0, 1 - 2r_n^{-1}\dist(x, \image(\beta_n)\cap B(z_n,r_n/2))\}
 \end{equation*}
 and
 \begin{equation*}
  \pi_n(x):= r_n^{-1}d(x,z'_n),
 \end{equation*}
 where $z'_n\in\image(\beta_n)$ lies between $\beta_n(0)$ and $z_n$ at distance $r_n/2$ from $z_n$. It follows that
 \begin{equation*}
  \int_0^1(f_n\circ\gamma_n)(s) (\pi_n\circ\gamma_n)'(s)ds\geq 1.
 \end{equation*}
 This concludes the proof of this direction. Now suppose $X$ is Gromov hyperbolic and let $\lambda>0$.  Fix $r>0$ and let  $\gamma:[0,1]\to X$ be a closed $\lambda r$-Lipschitz curve.
 It is not difficult to see that there exists $S\in\intcurr_2(L^\infty(X))$ such that $\bdry S=\gamma_{\#}\Lbrack\chi_{[0,1]}\Rbrack$ and $\mass{S}\leq C\lambda r$ for some 
 constant $C$ which does not depend on $r$ and $\lambda$. By the definition of $\|S\|$ we then have 
 \begin{equation*}
  \left|\int_0^1(f\circ\gamma)(s)(\pi\circ\gamma)'(s)ds\right|= |S(1, f,\pi)| \leq \lip(f)\lip(\pi)\mass{S}\leq C\lambda^3r^{-1}
 \end{equation*}
 for all $\lambda r^{-1}$-Lipschitz functions $f,\pi: X\to\R$.
 This completes the proof.
\end{proof}

We are ready for the proof of the main proposition of this section.

\begin{proof}[Proof of \propref{proposition:tree-bilip}]
 Assume $X$ is not Gromov hyperbolic. By \lemref{lemma:degenerate-loops} there exist $\lambda,\delta\in(0,\infty)$, a sequence $r_n\nearrow\infty$ and 
 $\lambda r_n$-Lipschitz maps $\gamma_n:[0,1]\to X$ with $\gamma_n(1)=\gamma_n(0)$ and such that the cycles $T_n\in\intcurr_1(X)$ defined by 
 $T_n:=\gamma_{n\#}\Lbrack\chi_{[0,1]}\Rbrack$ satisfy
 \begin{equation*}
  T_n(f_n,\pi_n) \geq \delta\quad\text{ for every $n\in\N$}
 \end{equation*}
 for suitable $\lambda r^{-1}_n$-Lipschitz functions $f_n,\pi_n:X\to\R$ with $f_n(\gamma_n(0))=\pi_n(\gamma_n(0))=0$.
 Let $X_\delta$ be a thickening of $X$ which admits a quadratic isoperimetric inequality for curves. By \lemref{lemma:almost-optimal-decomposition-one-cycles}, $X_\delta$ 
 also admits a quadratic isoperimetric inequality for $\intcurr_1(X_\delta)$.
 It then follows from 
 \cite[Lemma 3.4]{Wenger-GAFA} that there exists $S_n\in\intcurr_2(X_\delta)$ with $\bdry S_n=T_n$,
 \begin{equation*}
  \mass{S_n}\leq C\mass{T_n}^2\leq C \lambda^2r_n^2,
 \end{equation*}
 and such that the sequence of metric spaces $(\spt S_n, r_n^{-1}d_{X_\delta})$ is equi-compact and equi-bounded.
 In the above inequality, $C$ denotes the constant of the quadratic isoperimetric inequality for $\intcurr_1(X_\delta)$.
 Each of the following statements holds up to a subsequence.
 By Gromov's compactness theorem there exists a compact metric space $(Y,d_{Y})$ and isometric 
 embeddings $\psi_n: (Z_n, r_n^{-1}d_{X_\delta})\hookrightarrow Y$, where $Z_n:= \spt S_n\cup \gamma_n([0,1])$.
Furthermore, the compactness and closure theorem for currents \cite[Theorems 5.2 and 8.5]{Ambr-Kirch-curr} imply that $\psi_{n\#}S_n$ converges weakly to some $S\in\intcurr_2(Y)$.
Finally, $\psi_n(\spt S_n)$ converges to a compact subset $Z\subset Y$ with respect to the Hausdorff distance. 
Note that $(Z,d_Y)$ is an asymptotic set of $X$ and that furthermore $\spt S\subset Z$.
We now show that $S \not= 0$. Indeed, 
we can use McShane's extension theorem to first extend $f_n, \pi_n$ to $X_\delta$ and then to construct $\lambda$-Lipschitz functions $\tilde{f}_n, \tilde{\pi}_n: Y\to\R$ for which
$\tilde{f}_n\circ\psi_n = f_n$ and $\tilde{\pi}_n\circ\psi_n=\pi_n$. By Arzel\'a-Ascoli theorem the $\tilde{f}_n$ and $\tilde{\pi}_n$ converge
uniformly to $\lambda$-Lipschitz functions $\tilde{f},\tilde{\pi}: Y\to\R$.
Integration by parts finally yields
\begin{equation*}
 \begin{split}
   T_n(\tilde{f}&\circ\psi_n, \tilde{\pi}\circ\psi_n)\\
   &= T_n(f_n,\pi_n) - T_n((\tilde{\pi}-\tilde{\pi}_n)\circ\psi_n,\tilde{f}_n\circ\psi_n) 
    + T_n((\tilde{f}-\tilde{f}_n)\circ\psi_n,\tilde{\pi}\circ\psi_n)\\
   &\geq \delta - \lambda^2\|\tilde{\pi}-\tilde{\pi}_n\|_\infty - \lambda^2\|\tilde{f}-\tilde{f}_n\|_\infty
 \end{split}
\end{equation*}
and consequently
\begin{equation*}
 \bdry S(\tilde{f},\tilde{\pi})= \lim_{n\to\infty}(\psi_{n\#}T_n)(\tilde{f},\tilde{\pi})\geq\delta>0.
\end{equation*}
This shows that $S\not=0$. Consequently, by Theorem 4.5 of \cite{Ambr-Kirch-curr}, there exists a biLip\-schitz map $\varphi:K\subset\R^2\to Z$ with $\lm^2(K)>0$.
This yields a contradiction with the hypothesis and therefore concludes the proof.
\end{proof}

\section{Statement and proof of the main theorem}\label{section:isop-gromov}
The following is the main theorem of this article.

\bt\label{theorem:isop-general-case}
  Let $X$ be an admissible geodesic metric space and suppose there exist $\varepsilon>0$ and $s_0>0$ such that every Lipschitz loop $\gamma$ in $X$ with $\length(\gamma)\geq s_0$
  bounds an $S\in\intcurr_2(L^\infty(X))$ with
  \begin{equation}\label{equation:optimal-constant-currents}
   \mass{S}\leq \frac{1-\varepsilon}{4\pi} \length(\gamma)^2.
  \end{equation}
  Then $X$ is Gromov hyperbolic and, in particular, has a thickening which admits a linear isoperimetric inequality for curves.
\et

For the definition of `admissible' and conditions which imply admissibility see \secref{section:coarse-isoperimetric}.
Recall furthermore \lemref{lemma:Linfty-fillings} which asserts that for every metric space $Y$ isometrically containing $X$ and for every $S'\in\intcurr_2(Y)$ with boundary 
$\gamma$ there is an $S\in\intcurr_2(L^\infty(X))$ with boundary $\gamma$ and such that $\mass{S}\leq\mass{S'}$. In this sense, the existence of $S\in\intcurr_2(L^\infty(X))$ 
for which \eqref{equation:optimal-constant-currents} holds is the weakest condition we can ask. Note also that the theorem applies in particular to spaces in which loops in general
do not bound chains (such as Cayley graphs of groups).

\br\label{remark:change-of-measure}\rm
 Statements analogous to that in \thmref{theorem:isop-general-case} hold when mass $\mathbf{M}$ in \eqref{equation:optimal-constant-currents} is replaced
 by the parametrized Hausdorff or the Holmes-Thompson area, provided one works in the class of singular Lipschitz discs in $L^\infty(X)$ instead of integral currents. See the note after the proof.
\er

\begin{proof}[{Proof of \thmref{theorem:isop-general-case}}]
The proof is by contradiction and we therefore assume that $X$ is not Gromov hyperbolic. By \propref{proposition:tree-bilip} there exists a Lipschitz map
$\varphi:K\to (Z,d_Z)$ with $K\subset\R^2$ compact and $(Z,d_Z)$ an asymptotic subset of $X$ for which $\hm^2(\varphi(K))>0$. 
Let $\|\cdot\|$ be a norm on $\R^2$ as in \thmref{theorem:kirchheim-bilip-finite-set} and set $V:=(\R^2, \|\cdot\|)$. By approximation we may assume that the unit ball of $V$ is
the convex hull of finitely many points.
Let $\mathbb{I}_V\subset V$ be an isoperimetric subset of $V$ as in \secref{section:area-functionals}, 
set $a:=\length(\bdry\mathbb{I}_V)$ and let $\gamma:[0,a]\to \bdry\mathbb{I}_V$ be a 
parametrization by arc-length. Choose $M\in\N$ large enough (as below) and define
\begin{equation*}
 \Lambda:= \left\{\gamma(t_j): j=0,1,\dots, 2^M\right\},
\end{equation*}
where $t_j:= 2^{-M}aj$.
By the conclusion of \thmref{theorem:kirchheim-bilip-finite-set} and the definition of Gromov-Hausdorff limit there exists an $s_1\geq 10s_0$ arbitrary large and a 
$(1+\delta)$-biLipschitz map
$\psi:(\Lambda,\|\cdot\|)\to (X, \frac{1}{s_1}d)$. Here, we choose $\delta>0$ sufficiently small (see below). 
Let $X'$ denote the metric space $(X, \frac{1}{s_1}d)$ and note that by hypothesis, for every Lipschitz loop $c:[0,1]\to X'$ satisfying $\length(c)\geq\frac{s_0}{s_1}$ there exists
$S\in\intcurr_2(L^\infty(X'))$ with $\bdry S = c_{\#}\Lbrack\chi_{[0,1]}\Rbrack$ and
\begin{equation}\label{equation:rescaled-fillarea-ineq}
 \mass{S}\leq\frac{1-\varepsilon}{4\pi}\length(c)^2.
\end{equation}
Let now $c:[0,a]\to X'$ be 
a $(1+\delta)$-Lipschitz loop satisfying
\begin{equation*}
 c(t_j)=\psi(\gamma(t_j))\qquad\text{for all $j\in\{0,1,\dots,2^M\}$}
\end{equation*}
and let $T\in\intcurr_1(X')$ be given by $T:=c_\#\Lbrack\chi_{[0,a]}\Rbrack$.
Observe that 
\begin{equation*}
 \mass{T}\leq\length(c)\leq (1+\delta)a.
\end{equation*}
Since the unit ball of $V$ is the convex hull of finitely many points it follows from \cite[Lemma 9.19]{Thompson} that there exists an $n\in\N$ and a linear isometric embedding
$\varrho: V\hookrightarrow \ell^\infty_n$. Since $\ell^\infty_n$ is an injective metric space, $\varrho\circ\psi^{-1}$ can be extended to a $(1+\delta)$-Lipschitz map $\eta:L^\infty(X')\to\ell^\infty_n$. 
It is clear that for each $j=0, 1,\dots, 2^{M}-1$ there exists an $R_j\in\intcurr_2(\ell^\infty_n)$ satisfying
\begin{equation*}
 \bdry R_j= (\eta\circ c)_\#\Lbrack\chi_{[t_j,t_{j+1}]}\Rbrack - (\varrho\circ\gamma)_\#\Lbrack\chi_{[t_j,t_{j+1}]}\Rbrack
\end{equation*}
and
\begin{equation*}
 \mass{R_j}\leq C\left[2^{1-M}(1+\delta)^2a\right]^2,
\end{equation*}
where $C$ denotes the isoperimetric constant for $\intcurr_1(\ell^\infty_n)$. Set $R:=\sum_{j=0}^{2^M-1}R_j$ and let $S\in\intcurr_2(L^\infty(X'))$ 
be such that $\bdry S=T$ and $\mass{S}\leq C'\mass{T}^2\leq C'\length(c)^2$. Here, $C'$ is the isoperimetric constant for $\intcurr_1(L^\infty(X'))$.
Since $\bdry(\eta_\#S-R)=(\varrho\circ\gamma)_\#\Lbrack\chi_{[0,a]}\Rbrack$ and since $\mathbf M$ is semi-elliptic in the class of integral currents (see \secref{section:area-functionals})
it follows that
\begin{equation*}
 \mu^{m*}(\mathbb{I}_V)\leq \mass{\eta_\#S-R}\leq (1+\delta)^2\mass{S}+\mass{R}.
\end{equation*}
Since
\begin{equation*}
 \length(c)^2\leq (1+\delta)^2a^2=(1+\delta)^2\length(\bdry\mathbb{I}_V)^2\leq 4\pi(1+\delta)^2\mu^{m*}(\mathbb{I}_V)
\end{equation*}
this yields
\begin{equation}\label{equation:contradiction-fillarea}
 \length(c)^2\leq 4\pi(1+\delta)^4\mass{S} + 4\pi C(1+\delta)^62^{2-M}a^2. 
\end{equation}
Note that we also have
\begin{equation*}
 \length(c)\geq \frac{1}{\sqrt{C'}}\sqrt{\mass{S}}\geq \frac{s_0}{s_1}
\end{equation*}
if $s_1$ was chosen sufficiently large (only depending on $M$ and $\delta$).
From \eqref{equation:contradiction-fillarea} we conclude that
\begin{equation*}
 \mass{S}\geq \frac{1-\varepsilon/2}{4\pi}\length(c)^2
\end{equation*}
if $\delta$ is chosen small enough and $M$ large enough. Since this holds for all $S$ with $\bdry S=T$ this leads to a contradiction with 
\eqref{equation:rescaled-fillarea-ineq} and concludes the proof.
\end{proof}

By the facts stated in \secref{section:area-functionals} the above proof clearly works when mass $\mathbf{M}$ is replaced by $\mathcal{F}_{\mu^{ht}}$, provided one works in the
class of singular Lipschitz discs. Furthermore, since $\mathcal{F}_{\mu^{ht}}\leq\area$, the statement holds in particular for the parametrized Hausdorff area.

The following is a direct consequence of \thmref{theorem:isop-general-case} and gives a version of \thmref{theorem:isop-special-case} with the Hausdorff measure replaced by
the mass$*$-area.

\bc
 Let $X$ be a geodesic metric space and suppose there exists $\varepsilon>0$ such that every sufficiently long Lipschitz loop $\gamma$ in $X$
 bounds an $S\in\intcurr_2(X)$ with 
  \begin{equation}
    \mass{S}\leq \frac{1-\varepsilon}{4\pi}\length(\gamma)^2.
 \end{equation}
 Then $X$ is Gromov hyperbolic.
\ec

Note that we need not assume that $S$ is of disc type. Next we give the proof of \thmref{theorem:isop-special-case}.

\begin{proof}[Proof of \thmref{theorem:isop-special-case}]
 By \propref{proposition:admissible-thickening}, $X$ is admissible. 
 Now, the theorem follows directly from \thmref{theorem:isop-general-case} and \remref{remark:change-of-measure}.
\end{proof}

We end this section with the proof of the following theorem.

\bt\label{theorem:subquadratic-hyperbolic}
 Let $X$ be a geodesic metric space and suppose that for every $\nu>0$ there exists an $s_0>0$ such that every Lipschitz loop $\gamma$ in $X$ with 
 $\length(\gamma)\geq s_0$ bounds an $S\in\intcurr_2(L^\infty(X))$ satisfying
 \begin{equation*}
  \mass{S}\leq \nu\length(\gamma)^2.
 \end{equation*}
 Then $X$ is Gromov hyperbolic.
\et
It is important to note that we do not make the assumption that $X$ be admissible. Presently, it is not known to what extent the condition that $X$ be admissible can be
relaxed in \thmref{theorem:isop-general-case}. We mention that in general it is a difficult problem to determine the filling area in $L^\infty(X)$ of loops in $X$, see
e.g.\ \cite{Croke-Katz-filling-conj}.

\begin{proof}
 This follows directly from \lemref{lemma:degenerate-loops} together with the basic estimate \eqref{equation:H-basic-estimate}.
\end{proof}

\section{The sharp constant for the filling radius inequality}\label{Section:fillrad}
In this final section we determine the largest constant in a linear filling radius inequality in an admissible geodesic metric space which still implies Gromov hyperbolicity. 

Given metric spaces $X$ and $Y$, the filling radius in $Y$ of $T\in\intcurr_1(X)$ with $\bdry T=0$ is defined by
\begin{equation*}
 \fillrad_Y(T):= \inf\left\{r\geq 0: \text{$\exists S\in\intcurr_2(Y)$ with $\bdry S=T$ and $\spt S\subset B(\spt T, r)$}\right\}.
\end{equation*}
If $\gamma$ is a Lipschitz loop in $X$ then we write $\fillrad_Y(\gamma)$ for the filling radius in $Y$ of the integral current induced by $\gamma$.
The injectivity of $L^\infty(X)$ yields
\begin{equation*}
 \fillrad_{L^\infty(X)}(T)\leq \fillrad_Y(T)\leq \fillrad_X(T),
\end{equation*}
and in general these inequalities are strict. Indeed, if $\gamma$ parametrizes the unit circle in $\E^2$ then
$\fillrad_{\E^2}(\gamma) = 1$ and $\fillrad_{L^\infty(\E^2)}(\gamma)=\frac{\sqrt{3}}{2}$ as was shown by Katz in \cite{Katz}.

Next, let $\alpha_0$ be the largest number such that in any $2$-dimensional normed space $V$ there is a Lipschitz loop $\gamma: S^1\to V$ with $\length(\gamma)=1$ and
\begin{equation*}
 \fillrad_{L^\infty(V)}(\gamma)\geq \alpha_0.
\end{equation*}
It will be shown below that $\frac{3}{32}\leq \alpha_0\leq \frac{1}{8}$. We then have the following:

\bt\label{theorem:fillrad-hyperbolic}
 Let $X$ be an admissible geodesic metric space and suppose there exist $\varepsilon>0$ and $s_0>0$ such 
 that for every Lipschitz loop $\gamma$ in $X$ with $\length(\gamma)\geq s_0$
 \begin{equation*}
  \fillrad_{L^\infty(X)}(\gamma)\leq (1-\varepsilon)\alpha_0 \length(\gamma).
 \end{equation*}
 Then $X$ is Gromov hyperbolic and, in particular, has a thickening which admits a logarithmic filling radius inequality for curves.
\et

The theorem is clearly optimal in the class of admissible metric spaces, as follows from the definition of $\alpha_0$.
It generalizes results in \cite{Gromov-hyperbolic}, \cite{Drutu-french}, \cite{Papasoglu-fillrad} and improves the best known constant 
$\frac{1}{73}$ obtained by Papasoglu \cite{Papasoglu-fillrad}. The optimal value for the intrinsic filling radius inequality is conjectured to be $\frac{1}{8}$, see \cite{Papasoglu-fillrad}. 
At present we do not know the exact value of $\alpha_0$.

Before proving the theorem we show that $\frac{3}{32}\leq \alpha_0\leq \frac{1}{8}$. 
For this let $V$ be a normed space of dimension $k$ and recall that Jung's constant $J(V)$ is the smallest number $r\geq 1$ such that every set $A\subset V$ with $\diam A\leq 2$ is 
contained in some ball of radius at most $r$.
It is easy to see that $1\leq J(V)\leq 2$. Jung \cite{Jung} showed that $J(V)=1$ if and only if $V=\ell^\infty_k$. Bohnenblust \cite{Bohnenblust} furthermore proved 
that $J(V)\leq \frac{2k}{k+1}$. 
\bd 
 If $V$ is a $2$-dimensional normed space set
 \begin{equation*}
  \alpha_V:= \frac{1}{J(V)\length(\bdry B_V)},
 \end{equation*}
 where $B_V$ denotes the unit disc in $V$. Set moreover $\overline{\alpha}_0:= \inf_V\alpha_V$.
\ed
It is clear that for $V=\ell^\infty_2$ we have $\alpha_V=\frac{1}{8}$. 
The same holds for $\R^2$ endowed with the norm whose unit disc is a regular hexagon.
Go\l\c{a}b's Theorem \cite{Golab} asserts that $6\leq\length(\bdry B_V)\leq 8$ for every $2$-dimensional normed space $V$ 
and it thus follows from Bohnenblust's estimate that
\begin{equation*}
 \frac{3}{32}\leq \overline{\alpha}_0\leq \frac{1}{8}.
\end{equation*}
%

The estimate for $\alpha_0$ now follows directly from the proposition below.

\bp
 Let $V$ be a $2$-dimensional normed space and $\gamma:[0,1]\to V$ a Lipschitz parametrization of $\bdry B_V$. Then we have
 \begin{equation*}
  \fillrad_{L^\infty(V)}(\gamma)\geq \alpha_V\length(\gamma)
 \end{equation*}
 and consequently $\alpha_0\geq \overline{\alpha}_0$.
\ep

\begin{proof}
 Denote the norm on $V$ by $\|\cdot\|$ and abbreviate $B:=B_V$. 
 We first prove the proposition in the special case in which $B_V$ is the convex hull of points 
 \begin{equation*}
  \{\pm x_1,\dots,\pm x_n\}\subset V.
 \end{equation*}
 By \cite[Lemma 9.19]{Thompson} there then exists a linear
 isometric embedding $\varphi: V\to\ell^\infty_n$.
 Note that $\fillrad_{L^\infty(V)}(\gamma)= \fillrad_{\ell^\infty_n}(\varphi\circ\gamma)$ since
 $\ell^\infty_n$ and $L^\infty(V)$ are injective metric spaces. Set $T:=(\varphi\circ\gamma)_\#\Lbrack\chi_{[0,1]}\Rbrack$ and assume the existence of an $S\in\intcurr_2(\ell^\infty_n)$ with 
 $\bdry S=T$ and such that
 \begin{equation}\label{equation:spt-in-small-ball}
  \spt S\subset B(\spt T, r)\quad\text{ for some $r<\frac{1}{J(V)}$.}
 \end{equation}
 By the deformation theorem \cite[4.29]{Federer} and
 by the fact that $T$ is an integral polyhedral chain (see \cite[4.22]{Federer}) we may assume without loss of generality that $S$ is 
 an integral polyhedral chain satisfying \eqref{equation:spt-in-small-ball} and that each simplex in its support has diameter at most $\frac{1}{J(V)}-r$. 
 We now follow the arguments in the proof of Theorem 2 in \cite{Katz} in order to construct a Lipschitz retraction $\pi: \spt S\to \bdry B$ and
to arrive at a contradiction.
 Define $\pi$ on the $0$-skeleton of $\spt S$ by assigning to each vertex an arbitrary nearest point in $\bdry B$. 
 Consequently, if $A=\{u_1,u_2,u_3\}$ is the vertex set of a simplex then 
 \begin{equation*}
  \diam\pi(A)<\frac{2}{J(V)}
 \end{equation*}
 and hence $\pi(A)$ lies in a ball of radius strictly smaller than $1$ and thus in an open `hemisphere' of $\bdry B$. Therefore $\pi$ can be extended to a 
 Lipschitz map on the $1$-skeleton and then on all of $\spt S$ by sending the edges of simplices to the shortest paths on $\bdry B$ connecting the images of the 
 vertices and, furthermore, 
 `trivially' to the simplex. We conclude that $T=\pi_\# T= \bdry \pi_\# S=0$, which clearly contradicts the definition of $T$.

 As for the general case, choose $\lambda>1$ sufficiently close to $1$ and let $\{x_1,\dots,x_n\}\subset\bdry B$ be a finite subset such that the convex 
 hull $C$ of $\{\pm x_1,\dots,\pm x_n\}$ satisfies
 \begin{equation*}
  \frac{1}{\lambda}B\subset C \subset B.
 \end{equation*}
 Endow $\R^2$ with the norm whose unit ball is $C$ and denote this space by $W$. Clearly, the identity map $\psi: V \to W$ is $\lambda$-biLipschitz.
 We abbreviate $T:=\gamma_\#\Lbrack\chi_{[0,1]}\Rbrack$ and 
 $T':= \alpha_\#\Lbrack\chi_{[0,1]}\Rbrack$, where $\alpha:[0,1]\to W$ is a Lipschitz parametrization of $\bdry C$ with the same 
 orientation as that of $\gamma$, which we may assume to be counter-clockwise.
 By the special case already proved we have
 \begin{equation*}
  \fillrad_{L^\infty(W)}(T')\geq \frac{1}{J(W)}.
 \end{equation*}
 Let now be $S\in\intcurr_2(L^\infty(V))$ with $\bdry S=T$ and let $r\geq 0$ be such that 
 $\spt S\subset B(\spt T, r)$. If $\overline{\psi}: L^\infty(V)\to L^\infty(W)$ is a $\lambda$-Lipschitz extension of $\psi$ and if $\iota :\E^2\to W$
 denotes the identity map then $S':= \overline{\psi}_\#S - \iota_\#\Lbrack\chi_{B\ohne C}\Rbrack $
 satisfies $\bdry S'=T'$ and moreover
 \begin{equation*}
  \spt S'\subset B(\spt T', \lambda r+\lambda -1).
 \end{equation*}
 We conclude that
 \begin{equation*}
  \frac{1}{J(V)}\leq\frac{\lambda}{J(W)}\leq \lambda\fillrad_{L^\infty(W)}(T')\leq \lambda^2 r+\lambda(\lambda-1).
 \end{equation*}
 Since $\lambda>1$ was arbitrary this completes the proof.
\end{proof}

Finally, we give the proof of the above filling radius theorem. The strategy is analogous to that of \thmref{theorem:isop-general-case}. 

\begin{proof}[{Proof of \thmref{theorem:fillrad-hyperbolic}}]
Assume that $X$ is not Gromov hyperbolic and let 
$Z$, $\varphi$, $V$ and $\|\cdot\|$ be as in the proof of \thmref{theorem:isop-general-case}. 
Let $\gamma:[0,1]\to V$ be a Lipschitz loop of length $1$, parametrized by arc-length, for which
\begin{equation*}
 \fillrad_{L^\infty(V)}(\gamma)\geq \alpha_0.
\end{equation*}
Let $M\in\N$ be large enough and $\delta>0$ sufficiently small (as chosen below),
set $t_j:= \frac{j}{2^M}$ and
\begin{equation*}
 \Lambda:= \left\{\gamma(t_j): j=0,1,\dots, 2^M\right\}
\end{equation*}
and note that there exist just as in the proof of \thmref{theorem:isop-general-case} an $s_1\geq 10s_0$ arbitrary large and a $(1+\delta)$-biLipschitz map
\begin{equation*}
 \psi: (\Lambda, \|\cdot\|)\to X',
\end{equation*} 
where $X':=(X,\frac{1}{s_1}d)$.
Pick a $(1+\delta)$-Lipschitz loop $c:[0,1]\to X'$ satisfying
\begin{equation*}
 c(t_j)= \psi(\gamma(t_j))
\end{equation*}
for $j=0,1,\dots,2^M$ and note that 
\begin{equation*}
\frac{1}{1+2\delta}\leq \length(c)\leq (1+\delta)
\end{equation*}
if $M$ is large enough.
Let $\varphi:V\to\ell^\infty$ be a linear isometric embedding.
If $\eta:L^\infty(X')\to \ell^\infty$ is a $(1+\delta)$-Lipschitz extension of $\varphi\circ\psi^{-1}$ it is clear that
\begin{equation*}
 d(\varphi\circ\gamma(t),\eta\circ c(t))\leq 2^{-M}(1+\delta)^2
\end{equation*}
and hence there exists an $R\in\intcurr_2(\ell^\infty)$ such that $\bdry R=(\eta\circ c)_\#\Lbrack\chi_{[0,1]}\Rbrack-(\varphi\circ\gamma)_\#\Lbrack\chi_{[0,1]}\Rbrack$
and
\begin{equation*}
 \spt R\subset B\left(\varphi(\gamma([0,1])), 2^{-M}(1+\delta)^2\right).
\end{equation*}
If $S\in\intcurr_2(L^\infty(X'))$ satisfies $\bdry S=T:=c_\#\Lbrack\chi_{[0,1]}\Rbrack$ and if $r\geq 0$ is such that $\spt S\subset B(\spt T, r)$
then $S':=\eta_\#S-R\in\intcurr_2(\ell^\infty)$ satisfies $\bdry S'=(\varphi\circ\gamma)_\#\Lbrack\chi_{[0,1]}\Rbrack$ and we conclude
\begin{equation*}
 \frac{1}{1+\delta}\alpha_0\length(c)\leq \fillrad_{L^\infty(V)}(\gamma)\leq (1+\delta)r + 2^{-M}(1+\delta)^2.
\end{equation*}
Note furthermore that 
\begin{equation*}
 \length(c)\geq \frac{1}{1+2\delta}\geq \frac{s_0}{s_1}.
\end{equation*}
Choosing $\delta>0$ sufficiently small and $M\in\N$ large enough this leads to a contradiction with the assumption that
\begin{equation*}
 \fillrad_{L^\infty(X')}(c)\leq (1-\varepsilon)\alpha_0\length(c).
\end{equation*}
This concludes the proof.
\end{proof}

\end{document}